\documentclass[12pt]{article}

\usepackage{amssymb}

\renewcommand{\geq}{\geqslant}
\renewcommand{\leq}{\leqslant}
\renewcommand{\ge}{\geqslant}
\renewcommand{\le}{\leqslant}

\newtheorem{theorem}{Theorem}

\newtheorem{corollary}{Corollary}
\def\proof{\noindent {\bf Proof.}\,\,}

\def\qed{\hfill\vbox
  {\hrule\hbox{\vrule height1.3ex\hskip0.8ex\vrule}\hrule}\medskip}
\def\twid{\char"7E }

\def\eref#1{$(\ref{#1})$}

\def\tref#1{Theorem~$\ref{#1}$}

\def\capt#1{\vskip -0.7\baselineskip\caption{#1}\vskip 0.2\baselineskip}

\def\({\bigl(}  \def\){\bigr)}
\let\<=\langle  \let\>=\rangle

\def\Aut{\mathrm{Aut}}    \def\Is{\mathrm{Is}}

\def\Per{\mathrm{Per}}
\def\ppmod{\kern-0.7em\pmod}

\def\B{\mathcal{B}}

\def\M{\mathcal{M}}
\def\Red#1{R_{#1}}

\title{On the number of Latin squares}

\author{
Brendan D. McKay and Ian M. Wanless\\
\small Department of Computer Science\\[-0.5ex]
\small Australian National University\\[-0.5ex]
\small Canberra, ACT 0200, Australia\\
\small\tt $\{$bdm,imw$\}$@cs.anu.edu.au}

\date{}

\begin{document}

\maketitle

\begin{abstract}
We (1) determine the number of Latin rectangles with 11 columns and each
possible number of rows, including the Latin squares of order~11,
(2) answer some questions of Alter by showing that the number of reduced Latin
squares of order $n$ is divisible by $f!$ where $f$ is a particular
integer close to $\frac12n$, (3) provide a formula for the number of
Latin squares in terms of permanents of $(+1,-1)$-matrices, (4) find
the extremal values for the number of 1-factorisations of $k$-regular
bipartite graphs on $2n$ vertices whenever $1\leq k\leq n\leq11$, (5) show
that the proportion of Latin squares with a non-trivial symmetry group
tends quickly to zero as the order increases.
\end{abstract}

\setlength{\baselineskip}{3.0ex}

\section{Introduction}\label{intro}
For $1\le k\le n$, a $k\times n$ \emph{Latin rectangle\/} is a
$k\times n$ array $L=(\ell_{ij})$ with entries from $\{1,2,\ldots,n\}$
such that the entries in each row and in each column are distinct.  Of
course, $L$ is a \emph{Latin square\/} if $k=n$.  We say that $L$ is
\emph{reduced\/} if the first row is $(1,2,\ldots,n)$ and the first
column is $(1,2,\ldots,k)^T$.  If $\Red{k,n}$ is the number of reduced
$k\times n$ Latin rectangles then $L_{k,n}$, the total number of
$k\times n$ Latin rectangles, is $n!\,(n{-}1)!\,\Red{k,n}/(n{-}k)!$.
We will sometimes write $\Red{n,n}$ as~$\Red{n}$ and $L_{n,n}$ as $L_n$.

The determination of $\Red{k,n}$, especially in the case $k=n$, has been
a popular pursuit for a long time.  The number of reduced squares up
to order~5 was known to Euler~\cite{euler} and Cayley~\cite{cayley}.
McMahon~\cite{macmahon} used a different method to find the same
numbers, but obtained the wrong value for order~5.  The number of
reduced squares of order~6 was found by Frolov~\cite{frolov} and later
by Tarry~\cite{tarry}.  Frolov~\cite{frolov} also gave an incorrect
count of reduced squares of order~7.  Norton~\cite{norton} enumerated
the Latin squares of order~7 but incompletely; this was completed by
Sade~\cite{sade0} and Saxena~\cite{saxena2}.  The number of reduced
squares of order~8 was found by Wells~\cite{wells}, of order~9 by
Bammel and Rothstein~\cite{bammel}.

The value of $\Red{10}$ was found first in 1990 by the amateur
mathematician Eric Rogoyski working on his home computer and in the
following year by the present first author.  The resulting joint
paper~\cite{rogoyski} also presented the number of Latin rectangles
with up to 10 columns.  Before he died in 2002, Rogoyski worked for
several years on the squares of order~11 but the computing power
available to him was inadequate, despite his approach being sound.
Given the advance in computers since then, we can now complete the
computations moderately easily.

Several explicit formulas for general $n$ are in the literature
(\cite{shao}, for example).  Saxena~\cite{saxena2} succeeded in using
such a formula to compute~$\Red{7}$.
We will give another very simple formula in Section~\ref{formulas}.
At the time of writing, not even the asymptotic value of $\Red{n}$
is known.  In the case of rectangles, the best asymptotic result
is for $k=o(n^{6/7})$, by Godsil and McKay~\cite{godsil}.

\section{Terminology}

It can be useful to think of a Latin square of order $n$ as a set of
$n^2$ triples of the form (row, column, symbol).  For each Latin
square there are six {\em conjugate\/} squares obtained by uniformly
permuting the coordinates in each of its triples.  For example, the
transpose of $L$ is obtained by swapping the row and column
coordinates in each triple.

An {\em isotopism\/\/} of a Latin square $L$ is a permutation of its rows,
permutation of its columns and permutation of its symbols. The
resulting square is said to be {\em isotopic\/} to $L$ and the set of
all squares isotopic to $L$ is called an {\em isotopy class}. In the
special case when the same permutation is applied to the rows, columns
and symbols we say that the isotopism is an {\em isomorphism}. An
isotopism that maps $L$ to itself is called an {\em autotopism} of $L$ and
any autotopism that is an isomorphism is called an {\em
automorphism}. The {\em main class} of $L$ is the set of squares which
are isotopic to some conjugate of $L$. Latin squares belonging to the
same main class are said to be {\em paratopic\/} and a map which
combines an isotopism with conjugation is called a {\em paratopism}.
A paratopism which maps a Latin square to itself is called an {\em
autoparatopism\/} of the square.

The number of isomorphism classes, isotopy classes and main classes
has been determined by McKay, Meynert and Myrvold~\cite{myrvold}
for $n\le 10$. Our computation does not allow us to extract this
information for $n=11$. However, we do show in Section~\ref{s:nosym} that
$L_n/(6n!^3)$ provides an increasingly accurate estimate of the number
of main classes as $n$ grows.

\section{The Algorithm}\label{algo}

Our approach is essentially that introduced by Sade~\cite{sade0}, adapted
to the computer by Wells~\cite{wells,wellsbook}, and slightly improved
by Bammel and Rothstein~\cite{bammel}.  It was also used by McKay and
Rogoyski~\cite{rogoyski}.  Given a $k\times n$ Latin
rectangle~$L$, we can define a bipartite graph $B(L)$ with vertices
$C\cup S$, where $C=\{c_1,c_2,\ldots,c_n\}$ represents the columns
of~$L$ and $S=\{s_1,s_2,\ldots,s_n\}$ represents the symbols.  There
is an edge from $c_i$ to $s_j$ if and only if the symbol~$j$ appears
in column~$i$ of~$L$.  Thus $B(L)$ is regular of degree~$k$.  Clearly
$B(L)$ does not determine $L$ in general, since it does not record
the order of the symbols in each column. For us this is an
advantage, since it means there are many fewer graphs than there
are Latin rectangles.

Given a regular bipartite graph $B$ on $C\cup S$ of degree $k$, let $m(B)$
be its number of 1-factorizations, counted without regard to the
order of the factors.  Obviously $m(B)$ is an invariant of the
isomorphism class of $B$.  In speaking of isomorphisms and 
automorphisms of such bipartite graphs, we will admit the possibilities
that $C$ and $S$ are preserved setwise or that they are exchanged.
(More complex mixings of $C$ and $S$ would, in principle, be possible 
in the case of disconnected graphs, but we have chosen to 
disallow them.)  Using this convention, let $\Aut(B)$ be the
automorphism group of $B$ and let $\B(k,n)$ be a set consisting of
one representative of the isomorphism classes of bipartite graphs $B$
on $C\cup S$ of degree~$k$.

The theoretical basis of our approach is summarized in the following
theorem.  Parts~1 and~3 were proved in~\cite{rogoyski} and 
part~2 can be proved along similar lines.

\begin{theorem}\label{main}$ $\\
1. The number of reduced $k\times n$ Latin rectangles is given by
   $$\Red{k,n} = 2nk!(n{-}k)! \sum_{B\in\B(k,n)} m(B)|\Aut(B)|^{-1}.$$
2. The number of reduced Latin squares of order $n$ is given by
   $$\Red{n} = 2nk!(n{-}k)! \sum_{B\in\B(k,n)} m(B)m(\bar B)|\Aut(B)|^{-1},$$
   where $\bar B$ is the bipartite complement (the complement in $K_{n,n}$)
   of $B$ and $k$ is any integer in the range $0\leq k\leq n$.\\
3. Let $B\in\B(k,n)$ for $k\ge 1$.  Let $e$ be an arbitrary edge of $B$.
   Then
   $$m(B) = \sum_F m(B-F),$$
   where the sum is over all 1-factors $F$ of $B$ that include $e$.
\end{theorem}

For each $k=1,2,\dots,11$ in turn we found $m(B)$ for all $\B(k,11)$
using \tref{main}(3) and were then able to deduce $\Red{k,11}$ from
\tref{main}(1).  The number of graphs in $\B(k,11)$ is 1, 14, 4196,
2806508 and 78322916, for $k=1,\ldots,5$, respectively.  For $k\geq6$
the graphs in $\B(k,11)$ are the bipartite complements of those in
$\B(11-k,11)$. The main practical difficulty was the efficient
management of the fairly large amount of data.  Two implementations
were written in a way that made them independent in all substantial
aspects (except for their reliance on nauty \cite{nauty} to recognise the
isomorphism class of some graphs).  For example, they used different
edges $e$ in applying Theorem~\ref{main}(3), so that generally
different subgraphs were encountered.  The execution time of each
implementation was about 2 years (corrected to 1 GHz Pentium III), but
they would have completed in under 2 months if about 3 GB memory had
been available on the machines used.

We also ran the computations for $n\le 10$ and obtained the same
results as reported in~\cite{rogoyski}.  We repeat those results, and
include the new results, in Table~\ref{tab1}.  It is unlikely that
$\Red{12}$ will be computable by the same method for some time, since
the number of regular bipartite graphs of order 24 and degree 6 is
more than~$10^{11}$.

Note that our value of $\Red{11}$ agrees precisely with the numerical
estimate given in~\cite{rogoyski}, where estimates of $\Red{n}$ were given
for $11\leq n\leq15$.

\begin{table}[htp]
\begin{center}
\begin{tabular}{ccr|ccr}
 $n$&$k$&$\Red{k,n}$ & $n$&$k$&$\Red{k,n}$\\
\hline
1&1&1  & 9&1&1\\
\cline{1-3}
2&1&1  & &2&16687\\
 &2&1  & &3&1034\,43808\\
\cline{1-3}
3&1&1  & &4&20\,76245\,60256\\
 &2&1  & &5&11268\,16430\,83776\\
 &3&1  & &6&12\,95260\,54043\,81184\\
\cline{1-3}
4&1&1  & &7&224\,38296\,79166\,91456\\
 &2&3  & &8&377\,59757\,09642\,58816\\
 &3&4  & &9&377\,59757\,09642\,58816\\
\cline{4-6}
 &4&4     & 10&1&1\\
\cline{1-3}
5&1&1     &  &2&1\,48329\\
 &2&11    &  &3&81549\,99232\\
 &3&46    &  &4&14717\,45210\,59584\\ 
 &4&56    &  &5&746\,98838\,30762\,86464\\
 &5&56    &  &6&8\,70735\,40559\,10037\,09440\\
\cline{1-3}
6&1&1     &  &7&1771\,44296\,98305\,41859\,22560\\
 &2&53    &  &8&42920\,39421\,59185\,42730\,03520\\
 &3&1064  &  &9&75807\,21483\,16013\,28114\,89280\\
 &4&6552  &  &10&75807\,21483\,16013\,28114\,89280\\
\cline{4-6}
 &5&9408      & 11&1&1\\
 &6&9408      & &2&14\,68457\\
\cline{1-3}
7&1&1         & &3&79\,80304\,83328\\
 &2&309       & &4&143\,96888\,00784\,66048\\
 &3&35792     & &5&75\,33492\,32304\,79020\,93312\\
 &4&1293216   & &6&9\,62995\,52373\,29250\,51587\,78880\\
 &5&11270400  & &7&24012\,32164\,75173\,51550\,21735\,52640\\
  &6&16942080 & &8&86\,10820\,43577\,87266\,78085\,83437\,51680\\
  &7&16942080 & &9&2905\,99031\,00338\,82693\,11398\,90275\,94240\\
\cline{1-3}
8&1&1         & &10&5363\,93777\,32773\,71298\,11967\,35407\,71840\\
 &2&2119      & &11&5363\,93777\,32773\,71298\,11967\,35407\,71840\\
 &3&1673792   & \\
 &4&4209\,09504   & \\
 &5&27206\,658048   & \\
 &6&33\,53901\,89568   & \\
 &7&53\,52814\,01856   & \\
 &8&53\,52814\,01856   & \\
\end{tabular}
\end{center}
\capt{Reduced Latin rectangles}
\label{tab1}
\end{table}

\section{Some divisibility properties of $\Red{n}$}\label{divprops}

Despite obtaining the same value repeatedly for $\Red{11}$ by applying
Theorem~\ref{main}(2) for different $k$ in two independent
computations, we sought to check our answer further by determining its value
modulo some small prime powers.  By means of the algorithms described
in~\cite{myrvold}, we computed representatives $L$ of all the isotopy
classes of Latin squares of order~11 for which the order of the
autotopism group $\Is(L)$ is divisible by 5, 7, or~11.  The numbers of
such isotopy classes are listed in Table~\ref{tab2}.  Since the number
of reduced squares in the isotopy class of $L$ is $n\,n!/|\Is(L)|$,
these counts imply that $\Red{11}$ equals 8515 modulo 21175, in
agreement with our computations.

\begin{table}[ht]
\begin{center}
\begin{tabular}{c|c}
 $|\Is(L)|$ & isotopy classes \\
\hline
      5 &  55621 \\
      7 &   8065 \\
     10 &    359 \\
     11 &     24 \\
     14 &    160 \\
     20 &    102 \\
     21 &     45 \\
     22 &     12 \\
     55 &      6 \\
     60 &      3 \\
   1210 &      1
\end{tabular}
\end{center}
\capt{Isotopy classes with certain group sizes}
\label{tab2}
\end{table}

We also have the following simple divisibility properties.

\begin{theorem}\label{divis} For each integer $n\geq1$,\\
1. $\Red{2n+1}$ is divisible by\/ $\gcd(n!\,(n{-}1)!\,\Red{n},(n{+}1)!)$.\\
2. $\Red{2n}$ is divisible by\/ $n!\,$.
\end{theorem}

\proof Consider $\Red{2n+1}$ first.  We define an equivalence relation on
reduced Latin squares of order $2n+1$ such that each equivalence
class has size either $n!\,(n{-}1)!\,\Red{n}$ or $(n{+}1)!\,$.
Let $A$ be the leading principal minor of $L=(\ell_{ij})$ of order $n$.

If $A$ is a (reduced) Latin subsquare, then
the squares equivalent to $L$ are those obtainable by
possibly replacing $A$ by another reduced subsquare,
permuting the $n$ partial rows
$(\ell_{i,n+1},\ell_{i,n+2},\ldots,\ell_{i,2n+1})$ for $1\le i\le n$,
permuting the $n-1$ partial columns
$(\ell_{n+1,j},\ell_{n+2,j},\ldots,\ell_{2n+1,j})$ for $2\le j\le n$
then permuting columns $n+1,n+2,\ldots,2n+1$ to put the first row into
natural order.  These $n!\,(n{-}1)!\,\Red{n}$ operations are closed
under composition and give different reduced Latin squares, so each
equivalence class has size $n!\,(n{-}1)!\,\Red{n}$.

If $A$ is not a Latin subsquare, the squares equivalent to $L$ are
those obtainable by applying one of the $(n{+}1)!$ isomorphisms in
which the underlying permutation fixes each of the points
$1,2,\dots,n$.  No isomorphism of this form can be an automorphism of
a square in which $A$ is not a subsquare
(see~\cite[Theorem~1]{myrvold}).  Hence the squares obtained are
different and the equivalence class has $(n{+}1)!$ elements.

The case of $\Red{2n}$ is the same except the second argument
gives $n!$ instead of $(n{+}1)!\,$.
\qed

\begin{corollary}\label{divcor} \
If\/ $n=2p-1$ for some prime $p$, then $\Red{n}$ is divisible by
$\lfloor(n-1)/2\rfloor!\,$.  Otherwise, $\Red{n}$ is divisible by
$\lfloor(n+1)/2\rfloor!\,$.
\end{corollary}

\proof
This follows from Table~\ref{tab1} for $n\le8$.  
For $n\ge 9$, note that $m\mathbin{|} (m{-}2)!$ for $m>4$ unless
$m$ is prime.
\qed

Note that, for $n\ge 12$, the corollary gives the best divisor that can
be inferred from Table~\ref{tab1} and Theorem~\ref{divis}, except that
$\Red{13}$ is divisible by $7!$ and not merely by $6!$.

\medskip

Alter~\cite{alter} (see also Mullen~\cite{mullen}) asked whether an
increasing power of two divides $\Red{n}$ as $n$ increases and whether
$\Red{n}$ is divisible by 3 for all $n\geq6$.  Theorem~\ref{divis}
answers both these questions in the affirmative.  Indeed it shows
much more --- that for any integer $m>1$ the power of $m$ dividing $\Red{n}$
grows at least linearly in $n$. That is, for each $m$ there
exists $\lambda=\lambda(m)>0$ such that $\Red{n}$ is divisible by
$m^{\lfloor\lambda n\rfloor}$ for all $n$.

Alter also asked for the highest power of two dividing $\Red{n}$, and
here we must admit our ignorance. It seems from the evidence in
Table~\ref{primefact} that the power grows faster than linearly,
but we were unable to prove this.

\begin{table}[htp]
\begin{center}
$\begin{array}{c|c}
 n&\mbox{Prime factorisation of } \Red{n}\\
\hline
2& 1\\
3& 1\\
4& 2^2\\
5& 2^3\cdot7\\
6& 2^6\cdot3\cdot7^2\\
7& 2^{10}\cdot3\cdot5\cdot1103\\
8& 2^{17}\cdot3\cdot1361291\\
9& 2^{21}\cdot3^2\cdot5231\cdot3824477\\
10& 2^{28}\cdot3^2\cdot5\cdot31\cdot37\cdot547135293937\\
11& 2^{35}\cdot3^4\cdot5\cdot2801\cdot2206499\cdot62368028479\\
\end{array}$
\end{center}
\capt{Prime factorisations of $\Red{n}$ for $n\leq11$.}
\label{primefact}
\end{table}

\section{A formula for $\Red{n}$}\label{formulas}

The literature contains quite a few exact formulas for $\Red{n}$, but
none of them appear very efficient for explicit computation (though
Saxena~\cite{saxena2} managed to compute $\Red{7}$ using such a formula).

Perhaps the simplest formulas are those in~\cite{shao}, which relate
$\Red{n}$ to the permanents of all 0-1 matrices of order~$n$.  Here we
give one that is very similar but uses $\pm 1$ matrices instead.
Unlike the inclusion-exclusion proof of~\cite{shao}, we give a simple
analytic proof.

\begin{theorem}\label{exact}
Let $p(z)$ be any monic polynomial of degree $n$ and let $\M_n$ be the
family of all $n\times n$ matrices over $\{-1,+1\}$.  Then
$$L_n = 2^{-n^2} \sum_{X\in\M_n} p(\Per X)\,\pi(X),$$
where $\Per X$ is the permanent of~$X$ and $\pi(X)$ is the product of
the entries of~$X$.
\end{theorem}

\proof
If $X=(x_{ij})$ is an $n\times n$ matrix of indeterminates, then by definition
$\Per X = \sum_{\sigma\in S_n} T_\sigma$ where $S_n$ is the symmetric
group and $T_\sigma=x_{1\sigma(1)}x_{2\sigma(2)}\cdots x_{n\sigma(n)}$.
If the polynomial $p(\Per X)$ is expanded in terms of monomials in the
$x_{ij}$, then the only monomial involving every $x_{ij}$ comes from
products $T_{\sigma_1}T_{\sigma_2}\cdots T_{\sigma_n}$ where the
permutations $\sigma_1,\sigma_2,\ldots,\sigma_n$ are the rows of a
Latin square.  That is, the coefficient of the only monomial with each
$x_{ij}$ having odd degree is the number of Latin squares.  Multiplying by
$\pi(X)$ turns the required monomial into the only one that has even
degree in each $x_{ij}$.
Now summing over $X\in\M_n$ causes this monomial to be multiplied by
$|\M_n|=2^{n^2}$ while all the other monomials cancel out.
\qed

\section{Extremal graphs with respect to $m(B)$}

In our computations we learned the values of
$m(B)$ for each graph $B\in\B(k,n)$ for $n\le 11$.  In
Table~\ref{tabminmax} we record the maximum and minimum values, and
the number of graphs (in the column headed ``\#'') that achieve the
minimum.  The maximum is achieved uniquely in all cases.  Of course,
for $k\le 1$ the result is trivial and when $k\ge n-1$ the unique
graph has $m(B)=\Red{n}$, so we omit these cases.

In most cases, the graphs maximizing $m(B)$ are the same as those with
the maximum number of perfect matchings, as listed in~\cite{wanless}.
The only exceptions are as follows,
where the notation is that used in \cite{wanless}:
\begin{itemize}
\item For $n=7, k=5$ the graph maximising $m(B)$
is $\overline{2J_2\oplus D_3}$; 
\item For $n=9, k=6$ the graph maximising $m(B)$
is $\overline{3J_3}$; 
\item For $n=10, k=4$ the graph maximising $m(B)$
is $J_4\oplus\overline{3J_2}$; 
\item For $n=11, k=4$ the graph maximising $m(B)$
is $J_4\oplus\overline{J_3\oplus D_4}$. 
\end{itemize}
In the first of these cases the cited graph does, according
to \cite{wanless}, maximise the number
of perfect matchings, but does not do so uniquely.

\begin{table}[htp]
\begin{center}
\begin{tabular}{ccrcr}
 $n$&$k$&$\min m(B)$&\# & $\max m(B)$\\
\hline
4 & 2 & 1 & 1 & 2 \\[0.5ex]
5 & 2 & 1 & 1 & 2 \\
 & 3 & 4 & 1 & 6 \\[0.5ex]
6 & 2 & 1 & 1 & 4 \\
 & 3 & 8 & 4 & 24 \\
 & 4 & 168 & 1 & 224 \\[0.5ex]
7 & 2 & 1 & 1 & 4 \\
 & 3 & 8 & 3 & 48 \\
 & 4 & 456 & 2 & 576 \\
 & 5 & 54528 & 1 & 55296 \\[0.5ex]
8 & 2 & 1 & 1 & 8 \\
 & 3 & 16 & 18 & 96 \\
 & 4 & 1120 & 1 & 13824 \\
 & 5 & 3\,06432 & 1 & 4\,02432 \\
 & 6 & 2518\,94784 & 1 & 2583\,92064 \\[0.5ex]
9 & 2 & 1 & 1 & 8 \\
 & 3 & 16 & 7 & 288 \\
 & 4 & 2720 & 1 & 32256 \\
 & 5 & 17\,18784 & 1 & 23\,12192 \\
 & 6 & 35859\,25120 & 1 & 37975\,08096 \\
 & 7 & 2260\,68542\,91456 & 1 & 2271\,05054\,39232 \\[0.5ex]
10 & 2 & 1 & 1 & 16 \\
 & 3 & 24 & 2 & 576 \\
 & 4 & 6992 & 1 & 1\,29024 \\
 & 5 & 94\,57472 & 1 & 2167\,60320 \\
 & 6 & 4\,97127\,34208 & 1 & 7\,10221\,82400 \\
 & 7 & 92007\,32190\,63808 & 1 & 96252\,56413\,10208 \\
 & 8 & 51072\,82902\,02843\,87328 & 1 & 51411\,31576\,53646\,54080 \\[0.5ex]
11 & 2 & 1 & 1 & 16 \\
 & 3 & 32 & 25 & 1152 \\
 & 4 & 17040 & 1 & 3\,31776 \\
 & 5 & 494\,49728 & 1 & 15173\,22240 \\
 & 6 & 65\,69929\,07264 & 1 & 127\,45506\,81600 \\
 & 7 & 36\,18408\,76780\,25728 & 1 & 41\,31218\,87443\,35360 \\
 & 8 & 66\,74288\,35273\,45400\,70912 & 1 & 69\,04895\,67877\,90499\,02080 \\
 & 9 & 365\,09897\,56490\,71060\,26179\,78880 & 1 & 
           366\,51069\,03315\,59851\,95097\,12896 \\
\end{tabular}
\end{center}
\capt{Minimum and Maximum values of $m(B)$}
\label{tabminmax}
\end{table}

\section{Proportion of Latin squares with symmetry}\label{s:nosym}

In this section we prove that the proportion
of order $n$ Latin squares which have a non-trivial symmetry
tends very quickly to zero as $n\rightarrow\infty$.

\begin{theorem}\label{nosym}
The proportion of Latin squares of order $n$ which have a
non-trivial autoparatopy group is no more than
\begin{equation}\label{e:nosymprob}
n^{-3n^2/8+o(n^2)}.
\end{equation}
\end{theorem}

\proof 
Suppose that a Latin square $L=(\ell_{ij})$ of order $n$ has a non-trivial
autoparatopy group. Then by Lemma 4 in \cite{myrvold}, $L$ has a
autoparatopism $\alpha$ which fixes (pointwise) no more than one
quarter of the triples of $L$.

The number of possibilities for $\alpha$ is less than $6n!^3=o(n^{3n})$.
Given $\alpha$, we can construct each possible $L$ row by row.
Each entry is determined either by $\alpha$ and a previous entry,
or can be chosen in at most $n$ ways.  The latter possibility occurs
once per orbit of $\alpha$, and since $\alpha$ fixes at most $\frac14$
of the triples of $L$, the number of orbits is at most 
$(\frac14+\frac34\cdot\frac12)n^2=\frac58n^2$.
In total we find that there are most
$$o(n^{3n})n^{5n^2/8}$$
Latin squares with non-trivial autoparatopy group. Our result now 
follows immediately from the well known lower bound for $L_n$
(see, for example, Thm 17.2 in \cite{vLW})
that says that 
$$L_n\geq(n!)^{2n}n^{-n^2}\geq n^{n^2-o(n^2)}.$$
\qed

As a corollary to this last result, we see that the proportion of main
classes or of isotopy classes or of isomorphism classes whose members
have a non-trivial autoparatopy group is also bounded by
\eref{e:nosymprob}.  This is because each such class has somewhere
between $1$ and $6n!^3=o(n^{3n})$ members.

Another corollary is that the number of isomorphism classes, isotopy classes
and main classes of Latin squares
of order $n$ will be asymptotic to $L_n/n!$, $L_n/n!^3$ and $L_n/(6n!^3)$,
respectively.

  \let\oldthebibliography=\thebibliography
  \let\endoldthebibliography=\endthebibliography
  \renewenvironment{thebibliography}[1]{%
    \begin{oldthebibliography}{#1}%
      \setlength{\parskip}{0.4ex plus 0.1ex minus 0.1ex}%
      \setlength{\itemsep}{0.4ex plus 0.1ex minus 0.1ex}%
  }%
  {%
    \end{oldthebibliography}%
  }

\end{document}